\magnification=\magstep1
\input amstex
\documentstyle{amsppt}
\pageheight{8.5 truein}
\pagewidth{6.45 truein}
\NoBlackBoxes
\input epsf
\topmatter
\title
From factorizations of noncommutative polynomials to combinatorial
topology
\endtitle
\author
Vladimir Retakh
\endauthor
\rightheadtext{Noncommutative Polynomials and Topology}
\address
Department of Mathematics, Rutgers University,
Piscataway, NJ 08854
\endaddress
\email
vretakh\@math.rutgers.edu
\endemail
\dedicatory
To D.B. Fuchs on the occasion of his 70-th anniversary
\enddedicatory
\abstract
This is an extended version of a talk given by the author at the conference ``Algebra and Topology in Interaction"
on the occasion of the 70th Anniversary of D.B. Fuchs at UC Davis in September 2009.
It is a brief survey of an area originated around 1995 by I. Gelfand and the speaker.
\endabstract 

\endtopmatter
\document

\head 0. Introduction\endhead

Factorization of polynomials is one of the most fundamental and most classical
problems in mathematics. We know a lot about factorizations of polynomials over different fields
and, more generally, over commutative rings.
Much less is known about factorizations
of polynomials over noncommutative rings, e.g. polynomials with matrix
coefficients. Unlike their commutative counterparts, noncommutative
polynomials admit many different factorizations. This makes their theory
much harder and more interesting.

In 1995 I. Gelfand and the speaker constructed $n!$ different factorizations
of a noncommutative polynomial in one variable with $n$
roots in ``generic" position. Later with R. Wilson we studied ``algebras of pseudo-roots"
or ``noncommutative splitting algebras" associated with such factorizations.
These algebras can be described in terms of special directed graphs (quivers) called
layered graphs.

To any cell complex one can also associate a layered graphs and a ``splitting
algebra" defined by this graph. There are surprising connections between
properties of cell complexes and related splitting  algebras. Here I
construct a bridge between noncommutative algebra related to
factorizations of polynomials and combinatorial topology.

\head I. Factorizations of noncommutative polynomials\endhead

\noindent
{\bf 1.1}. Let us start with a secondary school problem. You are given a square polynomial
$$x^2+a_1x+a_2$$
over an associative (noncommutative) unital ring and its roots $x_1, x_2$.
How do you generalize the Viet\'e formulas $-a_1=x_1+x_2$ and  $a_2=x_1x_2$?

Let us assume that the difference $x_1-x_2$ is invertible and set
$$x_{1,2}=(x_2-x_1)x_2(x_1-x_2)^{-1},\ \  
x_{2,1}=(x_1-x_2)x_1(x_2-x_1)^{-1}.$$
\proclaim {Proposition 1.1}
$$-a_1=x_{1,2}+x_1=x_{2,1}+x_2$$
and 
$$a_2=x_{1,2}x_1=x_{2,1}x_2$$
(note the order of multipliers in the last formula!)
\endproclaim

\midinsert
\centerline{
\epsfxsize=2in
\epsfbox{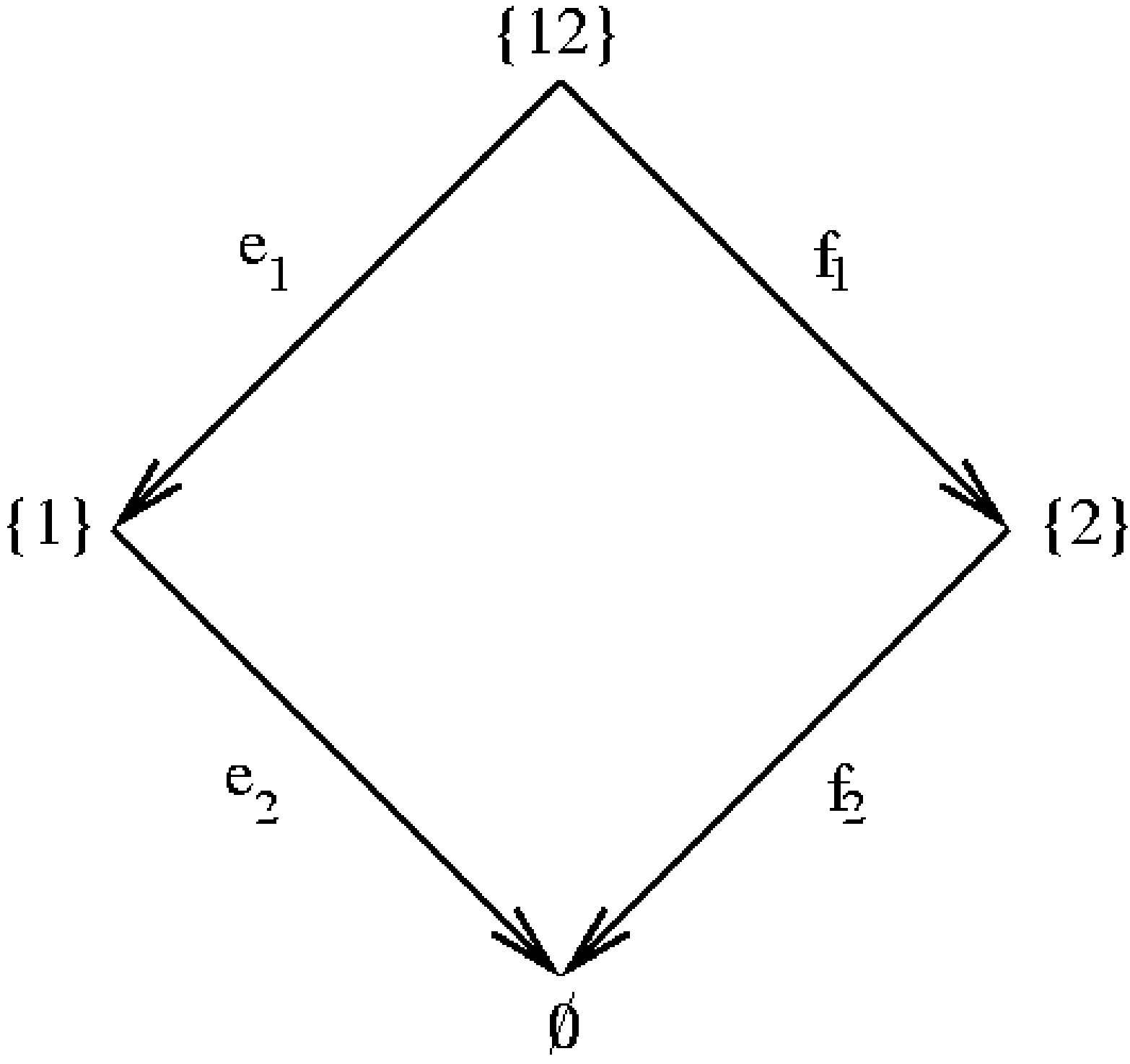}}
\botcaption{Figure 1}
\endcaption
\endinsert

Relations between $x_{1,2}, x_{2,1}, x_1, x_2$ can be described by the 
following graph (see Fig.~1). Here sums and products corresponding to the 
different paths from vertex $\{12\}$ to $\emptyset $ coincide. One can also present 
these relations in a unified way by setting $e_1=x_{1,2}$, $e_2=x_1$, $f_1=x_{2,1}$, $f_1=x_2$: 
$$(t-x_{1,2})(t-x_1)=(t-x_{2,1})(t-x_2)$$
where $t$ is a central variable.
\medskip \noindent
{\bf 1.2}. We proceed now to polynomials of any degree over a unital associative ring. Denote by $[1,n]$
the set of indices $\{1,2,\dots,n\}$. Assume that the polynomial
$$x^n+a_1x^{n-1}+\dots +a_{n-1}x+a_0$$
has $n$ roots $x_1,x_2,\dots ,x_n$ in a {\it generic} position, i.e. for any subset 
$\{i_1,i_2,\dots , i_{k+1}\}$ of $[1,n]$, $k\geq 1$, the Vandermonde matrix
$$W(i_1,i_2,\dots , i_{k+1})=\left ( \matrix
x_{i_1}^k&x_{i_2}^k&\dots &x_{i_{k+1}}^k\\
x_{i_1}^{k-1}&x_{i_2}^{k-1}&\dots &x_{i_{k+1}}^{k-1}\\
\dots &\dots &\dots &\dots \\
1 &1 &\dots &1 \endmatrix \right )$$ is invertible.

For the matrix $W(i_1,i_2,\dots , i_{k+1})$ denote the first row without
$x_{i_{k+1}}^k$ by $r$ and the last column without $x_{i_{k+1}}^k$ by $c$.
Set 
$$w_{i_1i_2\dots i_{k+1}}=x_{i_{k+1}}^k -r\cdot W(i_1,i_2,\dots , i_k)^{-1}\cdot c.$$

The element $w_{i_1i_2\dots i_{k+1}}$ is one of the {\it quasideterminants} of the matrix
$W(i_1,i_2,\dots ,i_{k+1})$ (see \cite {GR1-5, GGRW, RW1, Os1-2}). 

\medskip \noindent
{\bf Example 1.2}. One can see that $w_{i_1i_2}=x_{i_1}-x_{i_2}$.   

\medskip
The element $w_{i_1i_2\dots i_{k+1}}$ does not depend on the ordering of
the first $k$ indices $i_1,i_2,\dots i_k$, therefore we denote it by $w_{A_k, i_{k+1}}$
where $A_k=\{i_1, i_2,\dots , i_k\}$.

Assume now that the element $w_{A_k, i_{k+1}}$ is invertible and set
$$x_{A_k, i_{k+1}}=w_{A_k, i_{k+1}}\cdot x_{i_{k+1}}\cdot w_{A_k, i_{k+1}}^{-1}.$$

We are now ready to generalize Proposition 1.1.

\proclaim {Theorem 1.3} (\cite {GR4-5}) Let $i_1, i_2, \dots , i_n$ be an ordering of $1,2,\dots , n$.
Under our assumptions
$$-a_1=\sum _{k=1}^n x_{A_{k-1}, i_k},$$
$$a_2=\sum _{k>\ell} x_{A_{k-1}, i_k}x_{A_{\ell-1}, i_{\ell}}, \tag 1.1$$
$$\dots$$
$$(-1)^na_n=x_{A_{n-1}, i_n}x_{A_{n-2}, i_{n-1}}\dots x_{A_0, i_1}.$$
Here $A_0=\emptyset$ and $x_{A_0, i_1}=x_{i_1}.$
\endproclaim
\smallskip \noindent
{\bf Remark}.
In the commutative case the right-hand sides of these formulas are elementary
symmetric functions in $x_1,x_2,\dots , x_n$. In the noncommutative case the
right-hand sides are also symmetric in these variables (but you have to keep
the order of multipliers!). This remark is a cornerstone of the theory of 
noncommutative symmetric functions developed in \cite {GKLLRT} and subsequent papers. 
\smallskip

The formulas of Theorem 1.3 can be unified in the following way. Let $t$ be a central variable.
Set $P(t)=t^n+a_1t^{n-1}+\dots +a_{n-1}t+a_0$. Then
$$P(t)=(t-x_{A_{n-1},i_n})(t-x_{A_{n-2},i_{n-1}})\dots (t-x_{A_0,i_1})\tag 1.2 $$
and in a generic case there are $n!$ factorizations of $P(t)$ defined by
different orderings of $1,2,\dots ,n$. 

From now on we consider polynomials $P(t)$ over a ring $R$. An element $x\in R$
is a {\it left} root of $P(t)$ is $P(t)$ is divisible by $t-x$ from the left,
and a {\it right} root of $P(t)$ if $P(t)$ is divisible by $t-x$ from the right.
In this terminology every $x_{\emptyset, i}$ is a left root and every $x_{12\dots \hat i\dots n,i}$
is a right root of $P(t)$. 

Generally speaking, we cannot interpet elements $x_{A,i}$ as roots, and so we call them {\it pseudo-roots} due to
the following definition (see \cite {RSW3}).

\proclaim {Definition 1.4} An element $x\in R$ is a pseudo-root of a polynomial
$P(t)$ if there exist polynomials $S(t), T(t)\in R[t]$ such that
$$P(t)=S(t)(t-x)T(t).$$
\endproclaim
\smallskip \noindent
{\bf 1.3}. We construct now a universal algebra of pseudo-roots $Q_n$ defined by
generators and relations (see \cite {GRW1}). The algebra $Q_n$ is generated by
elements $x_{A,i}$ where
$A\subset [1,n]$ and $i\in [1,n]\setminus A$. The relations in
$Q_n$ are defined by the triples $(A,i,j)$ where $A\subset [1,n]$
and $i,j\in [1,n]\setminus A$, $i\neq j$:
$$x_{A\cup \{i\},j}+x_{A,i}=x_{A\cup \{j\},i}+x_{A,j},\tag 1.3a$$
$$x_{A\cup \{i\},j}\cdot x_{A,i}=x_{A\cup \{j\},i}+\cdot_{A,j}.\tag 1.3b$$

Relations (1.3) can be united by the equality
$$(t-x_{A\cup \{i\},j})(t-x_{A,i})=(t-x_{A\cup \{j\},i})(t-x_{A,j}) \tag 1.4$$
where $t$ is a formal central parameter.

Equality (1.4) implies that the products (1.2) where 
$A_k=\{i_1,i_2,\dots i_k\}$ do not depend on the orderings of $\{1,2,\dots , n\}$
and, therefore, define the monic polynomial $P(t)$ over $Q_n$. The generators
$x_{A,i}$ are the pseudo-roots of $P(t)$.

Note that relation (1.3) can be described geometrically by using
the diamond graph from Fig 1 where the sets $\{1,2\}$, $\{1\}$, $\{2\}$,
$\emptyset$ are replaced by $A\cup \{i,j\}$, $A\cup \{i\}$, $A\cup \{j\}$,
and $A$ respectfully. 

The next theorem shows that algebra $Q_n$ is a large algebra containing at least $n!+1$ 
free subalgebras with $n$ generators. 
Any ordered set $I=\{i_1,i_2,\dots , i_n\}$ defines a factorization (1.2) with
pseudo-roots $x_{A_{k-1},i_k}$, $k=1,2,\dots ,n$. Denote by $Q_I$ the subalgebra
generated by these pseudo-roots. 

\proclaim {Theorem 1.4} \cite {GRW1} For each $I$ algebra $Q_I$ is a free subalgebra of $Q_n$.
The intersection $Sym_n$ of free subalgebras $\bigcap Q_I$ over all orderings of $\{1,2,\dots , n\}$
is again the free subalgebra generated by elements $a_1, a_2,\dots , a_n$ defined by formulas
(1.1).
\endproclaim

One may consider elements $(-1)^ka_k$ as noncommutative analogs of the commutative elementary 
functions and algebra $Sym_n$ as an analogue of the corresponding commutative symmetric algebra. 
As we already said, a theory of noncommutative symmetric algebras was developed in a series
of papers starting with \cite {GKLLRT}.
\smallskip \noindent
{\bf 1.4}. To confirm that $Q_n$ has a large growth, we will pick a
set of linearly independent generators parametrized by subsets $B=(j_1,j_2,\dots , j_k)\subset [1,n]$.
The generators are $$x_B=x_{\emptyset, j_1}+x_{j_1,j_2}+\dots +x_{j_1,j_2,\dots j_{k-1}, j_k}.$$
Linear relations (1.3a) imply that $x_B$ does not depend on the order of indices $j_1,\dots , j_k$. 
Elements $x_B$ generate algebra $Q_n$ and satisfy only quadratic relations (see (\cite {GRW1}).
Therefore, $Q_n$ is a quadratic
algebra with a natural grading defined by $deg \; x_B=1$. Let $d_i$ be the dimension of the $i$-th
component of $Q_n$ according to this grading. Then according to \cite {GRW1, GGRSW} 
the graded dimension or the Hilbert series 
$H(Q_n, \tau)=\sum _{i\geq 0}d^i\tau ^i$ of the algebra $Q_n$ is

$$H(Q_n,\tau ) =\frac {1-\tau }{1 - \tau (2-\tau )^n}. \tag 1.5 $$

\head 2. From graphs to noncommutative polynomials\endhead

\noindent
{\bf 2.1}. In the previous section we used graphs to describe relations between pseudo-roots. 
Now we will use graphs to construct polynomials with noncommutative coefficients and their pseudo-roots
(see \cite {GRSW, GGRW, RSW3} and other papers).

Let $\Gamma =(V,E)$ be a directed graph (quiver) where $V$ is the set of vertices and $E$ is the set of edges.
For any $e\in E$ denote by $t(e)$ the tail of $e$ and by $h(e)$ the head of $e$.

Suppose that $\Gamma $ is a {\it layered} graph, i.e. $V=\coprod _{i=0}^n V_i$ and
$t(e)\in V_i$, $i>0$, implies $h(e)\in V_{i-1}$. We call $n$ the {\it height} of the graph
and we write $|v|=i$ if $v\in V_i$.
Set $V_+=\coprod _{i>0}V_i$.

\medskip \noindent
{\bf Examples 2.1}.  

\smallskip \noindent
a) Let $V$ be the set of subsets of the set $[1,n]$. For any $v\in V$
set $|v|$ to be the cardinality of $v$. An edge $e$ goes from $v$ to $w$ if $w$ is the subset of $v$
such that $|w|=|v|-1$. In this case $V_0=\{\emptyset\}$ and $V_n=\{[1,n]\}$.

\smallskip \noindent
b) Let $V$ be the set of linear subspaces of a vector space $F$ of dimension $n$.  
For any $v\in V$
set $|v|$ to be the dimension of $v$. An edge $e$ goes from $v$ to $w$ if $w$ is the subspace of $v$
such that $|w|=|v|-1$. In this case $V_0$ consists of the zero subspace and $V_n=\{F\}$.

\smallskip \noindent
c) This is a generalization of Example 2.1a. Let $X$ be a regular cell complex. Construct
the correspoinding layered graph $\Gamma _X=(V,E)$ as follows. Let $V$ be the set of cells $\sigma \in X$. 
Set $|\sigma|=\dim \/\sigma +1$.
An edge $e\in E$ goes from $\sigma$ to $\tau$ if and only if
$\tau$ lies in the closure of $\sigma$ and $|\tau|=|\sigma| -1$.
In this case $V_0={\emptyset}$ as in Example 1.    
\medskip \noindent
{\bf 2.2}. Following \cite {GRSW} we construct now an algebra $A(\Gamma)$ associated with a layered graph $\Gamma$. A path $\pi$ in $\Gamma$
is a sequence of edges $\pi=(e_1,e_2,\dots, e_k)$ such that $h(e_i)=t(e_{i+1})$ for $i=1,2,\dots, k-1$. We call
$t(e_1)$ the tail of $\pi$ and denote it by $t(\pi)$ and we call
$h(e_k)$ the head of $\pi$ and denote it by $h(\pi)$. The algebra $A(\Gamma)$ is generated over a given field by
generators $e\in E$ subject to the following relations. Let $t$ be a formal parameter commuting with edges $e\in E$. 
Any two paths $\pi=(e_1,e_2,\dots,e_k)$ and
$\pi'=(f_1,f_2,\dots,f_k)$ with the same tail and head define the relation
$$(t-e_1)(t-e_2)\dots (t-e_k)=(t-f_1)(t-f_2)\dots (t-f_k). \tag 2.1$$

In fact, relation (2.1) is equivalent to $k$ relations
$$e_1+e_2+\dots +e_k=f_1+f_2+\dots +f_k,$$
$$\sum _{i<j}e_ie_j=\sum _{i<j}f_if_j, \tag 2.2$$
$$\dots$$
$$e_1e_2\dots e_k=f_1f_2\dots f_k.$$
 
We call $A(\Gamma)$ the {\it splitting} algebra associated with graph $\Gamma$. The terminology is justified by
the following considerations. 
Assume that there are only one vertex $*$ of the minimal level $0$ and only one vertex $M$ of the maximal level $n$,
and that for any edge $e$ there exists
a path $\theta =(e_1,e_2,\dots ,e_n)$ from $M$ to $*$ containig $e$. Set 
$$P(t)=(t-e_1)(t-e_2)\dots (t-e_n).$$
 
Then $P(t)$ is a polynomial over $A(\Gamma)$, any edge in the graph $\Gamma$ corresponds to a pseudo-root
of $P(t)$, and any path from the maximal to minimal vertex corresponds to a factorization of $P(t)$ into a product of
linear factors. 
\smallskip \noindent
{\bf Example 2.2}. Let $\Gamma$ be the graph defined in Example 2.1a. Then $A(\Gamma)$ is 
the algebra $Q_n$ defined in Section 1.3.
\medskip \noindent
{\bf 2.3}. We will discuss now some properties of algebras $A(\Gamma)$. First we describe the
Hilbert series of such algebras
assuming that the layered graph $\Gamma$ has exactly one minimal vertex $*$ and
for any vertex $v\in V_i$, $i>0$  there is a path $\pi _v=(h_1,h_2,\dots, h_i)$ from $v$ to
the minimal vertex. The images of new elements $r_v=h_1+h_2+\dots +h_i$ in $A(\Gamma)$ 
do not depend on the choice of path from $v$ to $*$ and are linearly independent. 
According to formulas (2.2) elements $r_v$  generate algebra $A(\Gamma)$
and satisfy homogeneous relations of order two and higher. By assigning the degree one to
each element $r_v$ we obtain a  natural grading on $A(\Gamma)$. 
The corresponding Hilbert series $H(A(\Gamma),\tau )$ is well defined, provided that the set of vertices $V$ is finite.

To write an expression for the Hilbert series we introduce a {\it graded M\"obius function}
on $\Gamma$. We define a partial order on the set of vertices of $\Gamma$ by putting
$v>w$ if and only if there is a path from $v$ to $w$. The classical M\"obius function
$\mu$ on  $V\times V$ can be defined in the following way: $\mu(w,v)=0$ if $w\nless v$,
$\mu(v,v)=1$, and if $w<v$
$$\mu (w,v)=\sum_{w=v_0 < v_1 \dots < v_{\ell}=v} (-1)^{\ell}.$$
We introduce the {\it graded} M\"obius function of the graph as series
$$\Cal M(\Gamma , \tau)=\sum _{w<v}\mu (v,w)\tau ^{|v|-|w|}.$$

The following theorem was proved in \cite {RSW2}. 

\proclaim{Theorem 2.3} Let $\Gamma =(V,E)$ be a layered graph with a 
unique minimal element $*$ of level $0$. Assume that the set of vertices is finite.
Then

$$H(\Gamma, \tau ) =
\frac{1-\tau }{1-\tau M(\Gamma, \tau)}.$$
\endproclaim

One can show (see \cite {RSW4}) that the inverse series $H(\Gamma ,\tau)^{-1}$ is, in fact,
a polynomial. The degree of the polynomial equals to the height of $\Gamma$.
Theorem 2.3 implies formula (1.5) for algebra $Q_n$.

Generalized layered graphs corresponding to factorizations of polynomials into products of
not neccessary linear factors were studied in \cite {RW2}.

\head 3. Quadratic splitting algebras $A(\Gamma )$ and their dual algebras\endhead

\noindent
{\bf 3.1}. Let $\Gamma =(V,E)$ be a layered graph. As in Subsection 2.3 we assume that
the layered graph $\Gamma$ has exactly one minimal vertex and that
for any vertex $v\in V_i$, $i>0$  there is a path $\pi _v=(h_1,h_2,\dots, h_i)$ from $v$ to
the minimal vertex. In this case the splitting algebra $A(\Gamma)$ is defined by a set of homogeneous
relations of order $2$ and higher. Folowing \cite {RSW1} we describe now a condition when
algebra $A(\Gamma)$ is defined by quadratic relations only. 

We say that vertices $v,v'$ of the same level $i>0$ are connected by a {\it down-up}
sequence if there exist vertices $v=v_0,v_1,v_2,\dots ,v_k=v'\in V_i$ and 
$w_1,w_2,\dots, w_k\in V_{i-1}$ such that $w_j < v_{j-1}, v_j$ for $j=1,2,\dots , k$.
According to \cite {SRW1}, the layered graph $\Gamma $ is {\it uniform} if for any
pair of edges $e,e'\in E$ with a common tail, $t(e)=t(e')$, their heads $h(e), h(e')$
are connected by a down-up
sequence. It was proved in \cite {RSW1} that if the graph $\Gamma $ is uniform
then the splitting algebra $A(\Gamma)$ is quadratic, i.e. defined by relations of order $2$.

\medskip \noindent
{\bf 3.2}. For a quadratic algebra $A$ over a field $F$ there is a notion of the dual quadratic algebra
$A^!$. To define $A^!$ denote by $W$ the vector space of generators of $A$ and by $R\subset W\otimes W$
the linear space of relations of $A$. Denote by $W^*$ the dual space of $W$ and by
$R^{\perp}$ the annihilator of $R$ in $W^*\otimes W^*$. Algebra $A^!$ is a quadratic
algebra defined by generators from $W^*$ and relations from $R^{\perp}$. It is well-known (see, for example,
\cite {PP}) that an algebra $A$ is Koszul if and only if its dual algebra $A^!$ is. In this case their 
Hilbert series are connected by the simple equation 
$$H(A,\tau)H(A^!, -\tau)=1. \tag 3.1$$

We call a quadratic algebra $A$ {\it numerically Koszul} (see \cite {RSW4})
if the Hilbert series of algebras $A$ and $A^!$ satisfy equation (3.1).
\medskip \noindent
{\bf 3.3}. Assuming that a layered graph $\Gamma$ is uniform one can describe the dual algebra
$A(\Gamma )^!$ in terms of vertices and edges of the graph (see \cite {RSW3}). We describe
now a slightly different algebra $B(\Gamma )$.

There is a natural filtration on $A(\Gamma )$ defined by the ranking function $|\cdot |$.
The corresponding graded algebra is also quadratic.
Its dual algebra $B(\Gamma )$
can be described in the following way (see \cite {CPS}). For any $v\in V_+$ let $S(v)$ be the set of all vertices
$w\in V$ such that there is an edge going from $v$ to $w$. 

\proclaim {Theorem 3.1} The algebra $B(\Gamma )$ is generated by vertices $v\in V_+$ subject to the relations:

\noindent i) $u\cdot v=0$ if there is no edge going from $u$ to $v$;

\noindent ii) $v\cdot \sum _{w, w\in S(v)}w=0$.
\endproclaim

If the set of vertices $V$ is finite, the algebra $B(\Gamma )$ is finite-dimensional.
Algebras $B(\Gamma)$ were studied in \cite {CPS}.
According to the general theory $A(\Gamma )$ is Koszul if and only if $B(\Gamma)$ is Koszul.
If these algebras are Koszul then $H(B(\Gamma),\tau )=H(A(\Gamma), -\tau)^{-1}$.

Since not all splitting algebras are Koszul, we are going to investigate the difference
between $H(B(\Gamma),\tau )$ and $H(A(\Gamma), -\tau )^{-1}$ for layered graphs $\Gamma $.

Let $\Gamma =(V,E)$ be a layered graph, $V=\coprod _{i=0}^nV_i$. Denote by the same letter
$\Gamma$ the partially ordered set (poset) associated to the graph, by $\Delta (\Gamma)$ the {\it order}
complex of the poset [B], and by $H^i(\Gamma, F)$ the cohomology groups of $\Delta (\Gamma)$
with the coefficients in a field $F$. Let $b_i(\Gamma)$ be the dimension of 
$H^i(\Gamma, F)$. 

For a vertex $v\in V$ and $k\leq |v|$ denote by $\Gamma _{v,k}$ the graph induced by all vertices $w<v$,
$|w|\geq |v|-k+1$ with the added minimal vertex $*$, $|*|=0$.    

Note that $H(B(\Gamma ),\tau)$ is a polynomial in $\tau$ of degree $n$.

The following theorem was obtained in \cite {RSW5}.
\proclaim {Theorem 3.2} The coefficient of $\tau ^k$ in $H(A(\Gamma ),-\tau)^{-1}-H(B(\Gamma ),\tau)$ equals
$$\sum _{v\in V_{\ell}, \ell \geq k}\ \sum _{j=0}^{\ell-1}b_{\ell-j-1}(\Gamma _{v,k}).$$
\endproclaim  

\head 4. Splitting algebras associated to cell complexes\endhead

\noindent
{\bf 4.1}. Splitting algebras associated to cell complexes provide many
of examples of Koszul and non-Koszul algebras. The study of such algebras
was originated in \cite {RSW4} and continued in \cite {CPS, SS}.

Recall that to any regular cell complex $X$ one can associate the layered
graph $\Gamma _X$ (see Example 2.1c). 

\proclaim {Theorem 4.1} Let $X$ be a finite cell complex.
Algebras $A(\Gamma _X)$ and $B(\Gamma _X)$ defined over any ground field are Koszul.
\endproclaim

This theorem was proved in \cite {RSW4} for the abstract simplicial complexes
and in \cite {CPS} for any $CW$-complexes. The Koszulity of the algebra $Q_n$ was earlier proved in 
\cite {SW, Pi}.
Our construction gives a large number of Koszul algebras of exponential growth
as well as finite-dimensional Koszul algebras.
\medskip \noindent
{\bf 4.2}. Assume now that $X$ is a finite cell complex of dimension $n$.
We say that $X$ is {\it pure} if $X$ is the closure of its cells of dimension $n$.
Following \cite {CPS} we say that a pure complex $X$ is {\it connected by 
codimension one faces}
if any pair of cells of dimension $n$ is connected by a down-up sequence.  

Let $X$ be a pure cell complex of dimenion $n$ connected by codimension one faces.    
By adding to the graph $\Gamma _X$ one maximal vertex with the obvious
rank we obtain the graph $\Hat \Gamma _X$.

$A(\Gamma _X)$ and $A(\Hat \Gamma _X)$ have completely different properties.
The following theorem was proved
in \cite {RSW4} for two-dimensional manifolds and later in \cite {CPS} in general case.

\proclaim {Theorem 4.2} Let $X$ be a pure and regular $CW$-complex of dimension $n$
that is connected by codimension one faces. Let $F$ be a field.

\noindent
i) If algebras $A(\Hat \Gamma _X)$ and $B(\Hat \Gamma _X)$ over the ground field
$F$ are Koszul, then  $H^k(X,F)=0$ for all $0<k<n$. 

\noindent
ii) If $X$ is a manifold and $H^k(X,F)=0$ for all $0<k<n$ then
algebras $A(\Hat \Gamma _X)$ and $B(\Hat \Gamma _X)$ are Koszul over $F$.
\endproclaim

Theorem 4.2 implies that the Koszul property of algebras
$A(\Hat \Gamma _X)$ and $B(\Hat \Gamma _X)$ unlike the algebras
$A(\Gamma _X)$ and $B(\Gamma _X)$ over the ground field
depends on the field. As pointed in \cite {CPS}, if $X$ is a regular cellular representation of
$\Bbb {RP}^2$, then $A(\Hat \Gamma _X)$ and $B(\Hat \Gamma _X)$ are Koszul if
the field $F$ does not have characteristic $2$. 
 
\smallskip \noindent
{\bf Remark}.
Note that Theorem 3.2. can be used to describe the discrepancy between Hilbert series
of algebras $A(\Hat \Gamma _X, -\tau)^{-1}$ and $B(\Hat \Gamma _X, \tau)$ via Betti numbers
of the cell complex $X$. It shows that splitting algebras $A(\Hat \Gamma _X)$
and $B(\Hat \Gamma _X)$ carry a lot of information about topological properties
of the cell complex $X$. Note also that Theorem 3.2 can be used to prove a weaker version 
of Theorem 4.2.

\medskip \noindent
{\bf 4.3}. Theorem 4.2 gives a combinatorial solution to the question when
algebras $A\Hat \Gamma _X)$ and $B(\Hat \Gamma _X)$ are Koszul. These algebras 
are combinatorial invariants but not topological invariants. In a recent paper 
\cite {SS} it was shown, however, that the property of these algebras to be Koszul
is a topologigal invariant. The following theorem was proved in \cite {SS}.

\proclaim {Theorem 4.3} Let $X$ be a pure regular cell complex of dimension $n$  
connected through codimension one faces. Then algebras $A(\Hat \Gamma _X)$ and $B(\Hat \Gamma _X)$ are Koszul
for a field $F$ if and only if both of the following conditions hold:

\noindent
i) $\tilde H_i(X;F)=0$ for $i<n$;
\smallskip \noindent
ii) $\tilde H_i(X, X-\{p\};F)=0$ for each $p\in X$ and each $i<n$.
\endproclaim

Theorem 4.3 shows that the Koszul property for algebras $A(\Hat \Gamma _X)$ and $B(\Hat \Gamma _X)$
is a homeomorphic invariant. There are, however, homotopy equivalent, pure regular cell complexes $X$
and $Y$ such that $A(\Hat \Gamma _X)$ is Koszul but $A(\Hat \Gamma _Y)$ is not (see \cite {SS}).

\enddocument